\documentclass[11pt,a4paper]{article}

\usepackage[T1]{fontenc} 
\usepackage[utf8]{inputenc} 
\usepackage[frenchb,english]{babel} % le langage actif est english 
\usepackage{fourier}
\usepackage[scaled=0.875]{helvet}
\usepackage{courier}
\usepackage{dsfont}

\usepackage{amsmath,amsthm,amssymb,amsfonts,epic,latexsym,hyperref,natbib%,hypernat
}
\usepackage[dvips]{graphicx}
\usepackage[usenames]{color}

\newtheorem{theo}{Theorem}[section]
\newtheorem{defi}[theo]{Definition}
\newtheorem{prop}[theo]{Proposition}

\newtheorem{assu}{Assumption}

\newtheorem{exam}{Example}

\def\R{\mathbb{R}}
\def\N{\mathbb{N}}
\def\Z{\mathbb{Z}}

\def\w{\omega}
\def\P{\mathbb{P}} % proba environnement
\def\Pt{\mathbb{P}^\theta}
\def\Ps{\mathbb{P}^\star}
\def\E{\mathbb{E}} % esperance environnement
\def\Es{\mathbb{E}^\star}
\def\Et{\mathbb{E}^\theta}

\def\PP{\mathbf{P}} % proba annealed
\def\EE{\mathbf {E}} % esperance annealed
\newcommand{\PPt}{\mathbf{P}^{\theta}}

\newcommand{\PPs}{\mathbf{P}^{\star}}
\newcommand{\EEs}{\mathbf{E}^{\star}}

\def\1{\mathds{1}} % indicatrice

\def\t{\theta}
\newcommand{\ts}{{\theta^{\star}}}
\newcommand{\tz}{{\theta_{0}}}
\newcommand{\htet}{\widehat \theta_n}
\def\lsb{\ell^{\mathrm{sb}}}

\newcommand{\argmax}{\mathop{\rm Argmax}}

\def\dd{\mathrm{d}} %% pour le d des integrales
\def\ee{\mathrm{e}} %% pour le e exponentiel
\def\to{\rightarrow}

\def\Beta{\mathrm{B}}

\def\KL{d_{\mathrm{KL}}}
\def\nn{\nonumber}

\def\NB{\mathrm{NB}}
\def\NN{\mathcal{N}}

\def \eps{\varepsilon}

\def\to{\rightarrow}

\def\VV{\mathcal{V}}
\def\Beta{\mathrm{B}}
\def\proj{\mathrm{proj}}

\definecolor{lilas}{RGB}{182, 102, 210}
\begin{document}

\DeclareGraphicsExtensions{.pdf, .jpg, .jpeg, .png}

%%%%%%%%%%%%%%%%%%%%%%%%%%%%%
%%% TITRE, AUTEURS, RESUME, MOTS CLES %%%
%%%%%%%%%%%%%%%%%%%%%%%%%%%%%

\title{
 Maximum likelihood estimation in the context of a sub-ballistic random walk in a
 parametric random environment}

\author{Mikael  {\sc Falconnet}\footnote{Laboratoire  de Mathématiques
    et de Modélisation d'\'Evry, Universit\'e d'\'Evry Val
d'Essonne, UMR CNRS 8071, USC INRA,     E-mail:                          {\tt
  mikael.falconnet@genopole.cnrs.fr}; {\tt $\{$dasha.loukianova,arnaud.gloter$\}$@univ-evry.fr}
}
\and 
Arnaud {\sc Gloter}$^{*}$
\and
Dasha {\sc Loukianova}$^{*}$ 
}

 \maketitle

\begin {abstract}

We consider a one dimensional  sub-ballistic random walk evolving in a
parametric  i.i.d.\!   random environment.   We  study the  asymptotic
properties of the maximum  likelihood estimator (MLE) of the parameter
based on a  single observation of the path till the  time it reaches a
distant site.  In  that purpose, we adapt the  method developed in the
ballistic   case  by  \cite{Comets_etal}   and  \cite{FLM}.   Using  a
supplementary assumption  due to the specificity  of the sub-ballistic
regime, we  prove consistency and asymptotic normality  as the distant
site tends to infinity.  To emphazis 
% the difference between the two regimes and 
the role of the additional assumption, we investigate the Temkin model
with unknown support, and it turns out that the MLE is consistent but,
unlike in the ballistic regime, 
% the normalizedcriterion diverges and 
the  Fisher information  is infinite.  We also
explore the numerical performance of our estimation procedure.  
\end{abstract}

{\it Key words} :  Asymptotic normality, Sub-ballistic random walk, Confidence regions, Cramér-Rao efficiency, Maximum likelihood estimation, Random walk in random environment.
{\it MSC 2000} : Primary 62M05, 62F12; secondary 60J25.

%%%%%%%%%%%%%%%%%%%%
%%% INTRODUCTION %%%
%%%%%%%%%%%%%%%%%%%%

Let  $\w=(\w_x)_{x\in\Z}$  be a collection of independent  and  identically
distributed  (i.i.d.) $(0,1)$-valued random  variables
with  distribution  $\nu$. We suppose that the
law    $\nu=\nu_{\theta}$     depends    on    some  unknown  parameter
$\theta\in\Theta,$ where $\Theta \subset \R^d$ is assumed to be a compact set.
Denote  by  $\P^{\theta}=\nu_{\theta}^{\otimes   \Z}$  the  law  on
$(0,1)^{   \Z}$  of   the   environment  $\w$   and by 
$\E^{\theta}$ the expectation under this law.

For fixed  environment $\w$, let $X=(X_t)_{t\in\Z_+}$  be the Markov
chain on $\Z$ starting at $X_0=0$ and with transition probabilities 
\[
P_{\w}(X_{t+1}=y|X_t=x)=
\left \{
  \begin{array}{ll}
    \w_x&\mbox{if}\ y=x+1,\\
    1-\w_x&\mbox{if}\ y=x-1,\\
    0&\mbox{otherwise}.
  \end{array}
\right.
\]
The symbol $P_{\w}$ denotes the measure on the path space of $X$ given
$\w$, usually  called \emph{quenched} law. The  (unconditional) law of
$X$ is 
given by 
\[
  \PP^{\theta}(\cdot)=\int P_{\omega}(\cdot)\dd\P^{\theta}(\omega),
\]
this is the so-called
\emph{annealed}  law.   We  write  $E_{\w}$  and   $\EE^{\theta}$  for  the
corresponding  quenched and  annealed expectations,  respectively. The
behaviour of the process $X$ is related to the ratio sequence 
\begin{equation}
  \label{eq:rho}
  \rho_x=\frac{1-\w_x}{\w_x}, \qquad x \in {\mathbb Z},
\end{equation}
and we refer to~\cite{Sol} for the classification of $X$ between transient
or  recurrent cases  according  to whether  $\E^{\t}(\log \rho_0)$  is
different or not from $0$. 

The transient  case may  be further split  into two  sub-cases, called
\emph{ballistic} and \emph{sub-ballistic}  that correspond to a linear
and  a sub-linear speed  for the  walk, respectively.  More precisely,
letting $T_n$ 
be the first hitting time of the positive integer $n$, 
\begin{equation} \label{equa:HittingTime}
  T_n = \inf \{ t \in \N \, : \, X_t = n \},
\end{equation}
and   assuming  $\E^{\theta}(\log\rho_0)<0$   all   through,  we   can
distinguish the following cases. 
\begin{itemize}
\item[(a1)]  (Ballistic).  If   $\E^{\theta}(\rho_0)  <  1$,  then,  $
  \PP^{\theta} \mbox{-almost surely},$ 
%%  $T_n   /  n  \to  (1+\E^{\theta}(\rho_0))/(1-\E^{\theta}(\rho_0))$
%%  \mbox{$\PP^{\theta}$-a.s.}  
  \begin{equation}
    \label{eq:lln1}
    \frac{T_n}{n} \xrightarrow[n \to \infty]{} 
    \frac{1+\E^{\theta}(\rho_0)}{1-\E^{\theta}(\rho_0)}. 
  \end{equation}
\item[(a2)]  (Sub-ballistic). If  $\E^{\theta}(\rho_0)  \geq 1$,  then
  $T_n / n \to + \infty$, 
  $ \PP^{\theta} \mbox{-almost surely}$ when $n$ tends to infinity.
\end{itemize}

Moreover, the  fluctuations of $T_n$  depend in nature on  a parameter
$\kappa_\t \in (0,  \infty]$, which is defined as  the unique positive
solution of 
\begin{equation} 
  \label{equa:kappa}
  \E^{\theta}(\rho_0^{\kappa_\t})=1,
\end{equation}
when such a number exists, and $\kappa_\t=+\infty$ otherwise.  The 
sub-ballistic  case   corresponds  to  $\kappa_\t  \leq   1$.  In  our
statements, the quantity $\kappa_\t$ plays a crucial role that we will
emphazis when  it is implicitly  involved, since $\kappa_\t$  does not
appear explicitly in our assumptions.  

\cite{Comets_etal} provide a maximum likelihood estimator (MLE) of the
parameter of  the environment distribution  in the specific case  of a
transient \emph{ballistic} one-dimensional  nearest neighbour path. In
the  latter  work, the  authors  establish  the  consistency of  their
estimator while  the asymptotic  normality of the  MLE as well  as its
asymptotic  efficiency (namely,  that it  asymptotically  achieves the
Cramér-Rao bound)  is investigated in~\cite{FLM}.  The  method used in
these two  articles can  not be applied  directly for  a sub-ballistic
RWRE, due to the non-integrability  of the criterion function, but can
be  adapted  to the  sub-ballistic  regime.   However,  unlike in  the
ballistic regime,  the asymptotic behavior of the  estimator turns out
to be  very different when  estimating the support  of the law  of the
environment.  We   illustrate  this  when  we  consider  the
one-parameter Temkin model, a simple framework with finite and unknown
support,  which already reveals  the main  features of  the estimation
problem. One explanation  is that in the sub-ballistic  regime, due to
the  existence of  deeper local  traps of  the potential  than  in the
ballistic regime, the  walk spends a long time in  the bottom of these
traps,  and the  Fisher information  of the  support parameter
  becomes infinite. The non-finiteness of the Fisher
information suggests  that the convergence  of $\htet$ is  faster than
$\sqrt{n}$ and we provide  a simulation experiment that supports this.
Determining the true rate of convergence is a challenging problem that
we leave to further research.

This article is organised as follows.  
In  Section~\ref{sect:M_Estimator}, we  present our  MLE  procedure to
infer   the  parameter  of   the  environment   distribution  inspired
from~\citeauthor{Comets_etal}  and   recall  briefly some  already
known  results  on  an   underlying  branching  process  in  a  random
environment related to the RWRE. 
Then,   we  state  in   Section~\ref{sect:res}  our   consistency  and
asymptotic normality results, and present three examples of
environment    distributions     which    are    already    introduced
in~\cite{Comets_etal}  and~\cite{FLM}. The  MLE is  consistent  in the
three frameworks, but asymptotically  normal and efficient only in the
first  two cases.   In the  last  example, the  Fisher information  is
infinite and one of our assumptions fails.  
In Section~\ref{sect:proofs}, all the proofs are presented, and we
conclude with some simulation experiment in Section~\ref{sect:simus}.
%%

%%%%%%%%%%%%%%%%%%%%%%%%%%%%%%%
%%% CONSTRUCTION ESTIMATEUR %%%
%%%%%%%%%%%%%%%%%%%%%%%%%%%%%%%

\section{Maximum likelihood estimator  in the sub-ballistic transient case}
\label{sect:M_Estimator}

We always assume that $\Theta$ satisfies the following assumption.

\begin{assu} \label{assu:rho}
For any $ \theta\in\Theta$, 
\begin{itemize}
\item[i)] 
  $\E^{\theta}|\log\rho_0|<\infty$, 
  %%\quad
\item[ii)] \label{assu:trans}
  $\E^{\theta}(\log\rho_0)<0$,
\item[iii)] \label{assu:sousbal}
 $ \E^{\t}(\rho_0) \in [1,+\infty)$.
\end{itemize}
\end{assu}

%%In this section we briefly recall the estimator  proposed
%% in~\cite{Comets_etal} in the case  of the ballistic RWRE to infer a
%% parameter $\theta$ of the law of the environment. We explain why it
%%  doesn't work  anymore in  the sub-ballistic  case and  propose its
%% modification to recover this last case.  
 
The estimator  in~\cite{Comets_etal} is based  on the sequence  of the
number of  left steps performed by  the process $X$ from  sites $0$ to
site  $n$ at time  $T_n$  defined by~\eqref{equa:HittingTime}.
 More precisely, their estimator is the maximizer of the 
criterion function
\begin{equation}\label{eq:critbal}
\t \mapsto \ell_n(\t)=\sum_{x=0}^{n-1}  \phi_{\theta}(L_{x+1}^n,L_{x}^n),
\end{equation}
where $\phi_{\theta}$ is  the  function  from $\Z_+^2$ to $\R$ defined by
\begin{equation} \label{eq:PhiTheta}
   \phi_{\theta}(u,v) = \log\int_0^1a^{u+1}(1-a)^{v} \dd\nu_{\theta}(a),
\end{equation}
and  for any $x \in \{0,\ldots,n\}$
\begin{equation} \label{equa:leftstep}
  L_x^n := \sum_{t=0}^{T_n - 1} \1\{X_t=x, X_{t+1}=x-1 \}.
\end{equation}
\cite{Comets_etal} show that the limiting behavior of the
sequential log-likelihood function in  the case of ballistic RWRE is
equivalent to  \eqref{eq:critbal}. Recall from \cite{KKS}  that for an
i.i.d. environment, 
under  the annealed  law  $\PPt$, the  sequence $L_n^n$,  $L_{n-1}^n$,
$\dots$, $L^n_0$ has the same distribution as a branching process with
immigration  in  random  environment  (BPIRE)  denoted  $Z_0$,  $Z_1$,
$\dots$, $Z_n$ and defined by 
\begin{equation} \label{eq:BPRE}
Z_0=0, \quad \mbox{and for } k=0,\dots,n-1,\quad
Z_{k+1}=\sum_{i=0}^{Z_k}\xi_{k+1,i},
\end{equation}
with $\{\xi_{k,i}\}_{k\in\N;i\in\Z_+}$ independent and
\[
  \forall m\in\Z_+, \quad P_{\w}(\xi_{k,i}=m)=(1-\w_{k})^m \w_{k}.
\]
Under point~\hyperref[assu:trans]{$ii)$} of Assumption~\ref{assu:rho},
\citeauthor{Comets_etal} proved that the process $(Z_n)_{n\in\Z_+}$ is
a positive recurrent Markov chain with transition
kernel~$Q_\t$ defined as
\begin{equation} 
  \label{equa:Q}
  Q_{\t}(u,v)= \binom{u+v}{v} \int_0^1 a^{u+1}(1-a)^v \dd\nu_{\t}(a) =
  \binom{u+v}{v}\ee^{\phi_\t(v,v)}, \quad \forall u,v\in \Z_+.
\end{equation}
The unique invariant probability measure $\pi_\t$ of the process
$(Z_n)_{n\in\Z_+}$ is defined as 
\begin{equation} 
  \label{equa:pi}
  \pi_\t (u) = \Et [S(1-S)^u], \quad \forall u \in \Z_+,
\end{equation}
where
\begin{equation}
  \label{equa:S}
  S  =  \left(  \sum_{k=0}^\infty \prod_{i=1}^k  \rho_i  \right)^{-1}=
  (1+\rho_1  +  \rho_1  \rho_2  +  \cdots +  \rho_1  \ldots  \rho_k  +
  \cdots )^{-1} \in (0,1). 
\end{equation}
%%
%Using this scheme,  \citeauthor{Comets_etal} rewrite the likelihood in
%a simpler form and the MLE problem is reduced to an M-estimation 
 Due  to the equality  in law between $(L^n_n,\ldots  L^n_0)$ and
 $(Z_0,\ldots, Z_n),$ the MLE problem for RWRE is reduced to the one for 
the irreducible  positive recurrent  homogeneous Markov
chain  $(Z_n)_n$.  Thanks  to an  ergodic theorem  for  Markov chains,
\citeauthor{Comets_etal} 
proved  that in  the  \emph{ballistic} transient  case the  normalized
criterion  $\ell_n(\cdot)/n$ converges  in probability  to  a limiting
function  $\ell(\cdot)$  with   finite  values.  The  former  limiting
function identifies  the true value  of the parameter  and consistency
follows.     In    the     \emph{sub-ballistic}     transient    case,
\citeauthor{Comets_etal} prove 
that  the limiting function  $\ell(\cdot)$ still  exists but  might be
infinite everywhere, and  hence do not identify the  true value of the
parameter.  
%%The former  fact is due  to the  ergodic theorem
%%which in the \emph{sub-ballistic} case might provide an infinite limit.  

%To modify properly the criterion function so that the ergodic theorem for
%Markov chains provides  a limiting function with finite  values in the
%\emph{sub-ballistic} case, 
Let us explain briefly where is the problem. 
Introduce  the probability  measure  $\tilde \pi_\t$  on $\Z_+  \times
\Z_+$ defined as 
\begin{equation} \label{equa:pi_tilde}
  \tilde  \pi_\t(u,v) = \pi_\t  (u) Q_\t(u,v),
\end{equation}
and  denote $\tilde \pi_{\t}(g)$  for any  function $g:\Z_+^2  \to \R$
such  that $  \sum_{x,y}\tilde \pi_\t  (x,y) |g(x,y)|  <  \infty$ with
$\tilde \pi_\t$ the quantity defined as 
\begin{equation}\label{eq:pi_tilde_g}
\tilde \pi_\t(g) = \sum_{(x,y)\in \N^2}\tilde \pi_\t (x,y) g(x,y).
\end{equation}
In \cite{Comets_etal}, the limiting function $\ell(\cdot)$ is defined
as $\t \mapsto \tilde \pi_{\ts}(\phi_\t)$ where $\ts$ is the true parameter
value,  and the  integrability of  $\phi_\t$ with  respect  to $\tilde
\pi_{\ts}$  is equivalent  to  the  existence of  a  first moment  for
$\pi_{\ts}$.   We will  see  in Proposition~\ref{prop:pimoments}  that
$\kappa_\t$ defined by~\eqref{equa:kappa}  is the upper critical value
for the existence of finite  moments for $\pi_\t$. Therefore, since in
the \emph{sub-ballistic} case, we  have $\kappa_{\ts} \leq 1$, we know
that  $\pi_{\ts}$ does  not  have  a first  moment  and $\ell(\t)$  is
infinite. 
In the light  of this, the natural idea is  to consider the difference
of two log-likelihood functions.

\begin{defi}
  Fix $\tz \in \Theta$. The criterium function
  $\t \mapsto \lsb_n(\t)$ is defined as 
  \begin{equation}\label{eq:l} 
    \lsb_n(\theta)=\sum_{x=0}^{n-1}\left[ \phi_{\theta}(L_{x+1}^n,L_{x}^n)-
      \phi_{\theta_0}(L_{x+1}^n,L_{x}^n)\right ].  
  \end{equation}
  An estimator $ \htet$ of $\theta$ is defined as a measurable choice
  \begin{equation}\label{eq:estimator}
    \htet \in \argmax_{\theta\in\Theta}\lsb_n(\theta).
  \end{equation}
\end{defi}
As soon as the function $\t \mapsto \phi_\t(u,v)$ is continuous on the
compact parameter set  $\Theta$ for any pair of  integers $(u,v)$, the
criterion  function  $\lsb_n(\cdot)$  achieves  its maximum,  and  the
estimator  $\htet$   is  well  defined   as  one  maximizer   of  this
criterion. However, it is not necessarily unique.  

%%%%%%%%%%%%%%%%%
%%% RESULTATS %%%
%%%%%%%%%%%%%%%%%
\section{Consistency and asymptotic normality results} \label{sect:res}

From now  on, we assume  that the process  $X$ is generated  under the
true parameter value  $\ts$, an interior point of  the parameter space
$\Theta$, that we  aim at estimating. We shorten  to $\PPs$ and $\EEs$
(resp.   $\Ps$  and  $\Es$)  the   annealed  (resp.  the  law  of  the
environment   )  probability   $\PP^{\ts}$   (resp.  $\P^{\ts}$)   and
corresponding expectation~$\EE^{\ts}$ (resp. $\E^\ts$) under parameter
value $\ts$.

\subsection{Consistency result}

Assumption~\ref{assu:consistance} below ensures  that the maximizer of
criterion $\lsb_n$ is a consistent estimator of the unknown parameter.  

%%%%%%%%%%%%%%%%%%%%%%%%%%%%%%%
%%% CONSISTENCY ASSUMPTIONS %%%
%%%%%%%%%%%%%%%%%%%%%%%%%%%%%%%

\begin{assu} \label{assu:consistance}
\ 

\begin{itemize}
\item[i)]\label{assu:cont} (Continuity). 
For  any $(x,y)\in  \N^2$, the  map $\t  \mapsto
\phi_{\theta}(x,y)$ is continuous on the parameter set $\Theta$.
 \item[ii)] \label{assu:ident} (Identifiability). For any 
 $ (\theta,\theta')\in\Theta^2,$ 
 $\nu_{\theta}\neq\nu_{\theta'}\iff \theta\neq\theta'.$
 \item[iii)]\label{assu:integrability} (Uniform integrability).  For 
   any $\t \in \Theta$, 
$\tilde  \pi_\t  \big(  \sup_{\theta'   \in  \Theta}  |  \phi_{\t'}  -
\phi_{\tz} | \big) < \infty$. 
\end{itemize}
\end{assu}

We now state our main result.

\begin{theo}(Consistency). \label{theo:consistency}
Under Assumptions~\ref{assu:rho} and~\ref{assu:consistance}, 
for any  choice of $\widehat\theta_n$ satisfying~\eqref{eq:estimator},
we 
have
\[
  \lim_{n\to\infty}\widehat\theta_n=\ts,
\]
in $\PPs$-probability.
\end{theo}
Theorem~\ref{theo:consistency} is a straight application of Theorem
5.7 in~\cite{VdV}. Hence, it suffices to check that the assumptions of
the former  theorem are fulfilled. The  first one is  the uniform weak
law of large numbers for the renormalized criterion given in
Proposition~\ref{prop:LimitEllTheta}, and the second one is the
statement       of      Proposition~\ref{prop:IdentificationEllTheta}.
Sections~\ref{sect:LimitEllTheta} and~\ref{sect:IdentificationEllTheta} 
are dedicated to their respective proof.  
\begin{prop} \label{prop:LimitEllTheta}
Under   Assumptions~\ref{assu:rho}   and~\ref{assu:consistance},   the
following uniform convergence holds: 
\begin{equation} \label{equa:ULLN}
\sup_{\t \in \Theta}  \left | \frac 1n \lsb_n(\t)  - \lsb(\t) \right |
\xrightarrow[n \to \infty]{} 0 \quad \mbox{in $\PPs$-probability}, 
\end{equation}
with
\begin{equation} \label{eq:defell}
  \lsb(\t) = \tilde \pi_\ts (\phi_{\t} - \phi_{\tz}).
\end{equation}
\end{prop}

\begin{prop} \label{prop:IdentificationEllTheta}
Under  Assumptions~\ref{assu:rho} and~\ref{assu:consistance}, for any
$\eps >0$, 
\begin{equation} \label{equa:ident}
  \sup_{\theta: \|\theta-\ts\| \geq \eps} \lsb(\theta)<\lsb(\ts).
\end{equation}
\end{prop}

From Section~\ref{sect:M_Estimator}, point~\hyperref[assu:trans]{$iii)$} of 
Assumption~\ref{assu:consistance}   is   essential   to  ensure   that
$\lsb(\cdot)$    takes    finite    values   and    therefore    prove
consistency.  This
point  can be  expressed in  terms of  the growth  of $\dot  \phi$ and
thereby  is related  to the  existence of  moments of  the probability
distribution  $\pi_{\ts}$  which   are  characterized  in  Proposition
\ref{prop:pimoments} below.  Note that in  the ballistic regime,
since $\pi_{\ts}$  possesses a finite  first moment and the  growth of
$\dot \phi$ is linear, point~\hyperref[assu:trans]{$iii)$} of 
Assumption~\ref{assu:consistance} is automatically satisfied. 

\begin{prop}\label{prop:pimoments} 
Let   $\kappa_\t$  defined   by~\eqref{equa:kappa}  and   $\alpha  \in
(0,+\infty)$.      Under     point~\hyperref[assu:trans]{$ii)$}     of
Assumption~\ref{assu:rho}, the following dichotomy holds: 
\begin{itemize}
\item[i)] \label{prop:pimomentfini}
 $\alpha < \kappa_\t \Longrightarrow \quad
 \sum_{k=0}^{\infty}k^{\alpha}\pi(k)<\infty$; 
\item[ii)] \label{prop:pimomentinfini}  
$\alpha \geq \kappa_\t \Longrightarrow\quad
\sum_{k=0}^{\infty}k^{\alpha}\pi(k)=\infty$. 
\end{itemize}
\end{prop}
 Section~\ref{sect:pimoments} is dedicated to the proof of
\ref{prop:pimoments}.

%%%%%%%%%%%%%%%%%%%%%%%%%%%%%%
%%% NORMALITE ASYMPTOTIQUE %%%
%%%%%%%%%%%%%%%%%%%%%%%%%%%%%%

\subsection{Asymptotic normality results} 

The asymptotic  normality result  in \cite{FLM} involves  the gradient
and  the   second  derivative  of  $\ell_n(\cdot)$   with  respect  to
$\t$. Since they  are equal to the gradient  and the second derivative
of $\lsb_n(\cdot)$ with respect to  $\t$, their result can be extended
to the sub-ballistic  case under the same assumptions  and without any
modification of their proof.  

In the following,  for any function $g_\t$ depending  on the parameter
$\t$, the symbols $\dot  g_\t$ or $\partial_\t g_\t$ and $\ddot{g}_\t$
or  $\partial_\t^2  g_\t$  denote  the (column)  gradient  vector  and
Hessian   matrix  with  respect   to  $\t$,   respectively.  Moreover,
$Y^\intercal$  is the row  vector obtained  by transposing  the column
vector $Y$.

\begin{assu} \label{assu:normality}
\

\begin{itemize} 
\item[i)] (differentiability). \label{assu:diff}
The collection of probability measures $\{ \nu_\t \, : \, \t \in\Theta
\}$ is such that for any $(x,y)\in \N^2$, the map $\t \mapsto 
\phi_{\theta}(x,y)$ is twice continuously differentiable on $\Theta$.
\item[ii)] (Regularity conditions).\label{assu:smoothness}
For any $\t \in \Theta$, there exists some $q>1$ such that
$\tilde \pi_\t\Big(\|\dot \phi_{\t}\|^{2q}\Big)<+\infty$. 
\item[iii)] (Invertibility). \label{assu:invertibility} 
For any $u \in  \Z_+$, $\sum_{v \in \Z_+} \dot Q_{\t}(u,v)=\partial_\t
\big( \sum_{v\in \Z_+} Q_{\t}(u,v)$ \big).  
%%
%%For any $x\in\N,$
%%\begin{equation}\label{eq:chaing_ddot}
%%\quad \sum_{y\in \N}\ddot Q_{\t}(x,y)=0.
%%\end{equation}
%%
\item[iv)] (Uniform conditions).\label{assu:uniform}
For   any   $\t   \in   \Theta$,  there   exists   some   neighborhood
$\mathcal{V}(\t)$ of $\t$ such that 
$
 \tilde   \pi_\t  \Big(\sup_{\t'   \in   \mathcal{V}(\t)}  \|\dot
 \phi_{\t'} \|^2 \Big) < +\infty$ 
and 
$  \tilde  \pi_\t \Big( \sup_{\t' \in 
   \mathcal{V}(\t)} \|\ddot \phi_{\t'} \|\Big) < +\infty$.
\item[v)] (Fisher information matrix).\label{assu:nonsingular}
  For  any  value  $\t  \in \Theta$,  the  matrix  $\Sigma_{\t}=\tilde
  \pi_{\t} \Big(\dot \phi_{\t}^{\phantom{\intercal}} \dot 
\phi_{\t}^\intercal \Big)= - \tilde \pi_\t (\ddot \phi_{\t})$ is  non
  singular.  
\end{itemize}
\end{assu}

\begin{theo}\label{thm:CLT_dotphi}
Under Assumptions~\ref{assu:rho} to~\ref{assu:smoothness}, the 
score vector  sequence $ {\dot  \ell}_n^{\mathrm{sb}}(\ts) / \sqrt{n}$
is asymptotically 
normal with mean zero and finite covariance matrix $ \Sigma_{\ts}$.
\end{theo}

\begin{theo} (Asymptotic normality). \label{theo:CLT}
Under   Assumptions~\ref{assu:rho}  to~\ref{assu:nonsingular},
for any  choice of $\widehat\theta_n$ satisfying~\eqref{eq:estimator},
the  sequence  $\{\sqrt{n}  (\htet-\ts)\}_{n  \in  \N}$  converges  in
$\PPs$-distribution  to   a  centered  Gaussian   random  vector  with
covariance  matrix  $   \Sigma_{\ts}  ^{-1}$.   
\end{theo}

Note  that  the limiting  covariance  matrix  of  $\sqrt{n} \htet$  is
exactly the inverse Fisher information matrix of the model.  As such, our
estimator is  efficient.

%%%%%%%%%%%%%%
%% EXEMPLES %%
%%%%%%%%%%%%%%

\subsection{Examples}
We  illustrate  our results  in  the  same  frameworks than  the  ones
presented by~\cite{Comets_etal}
and \cite{FLM}. Note that point~\hyperref[assu:integrability]{$iii)$} of 
Assumption~\ref{assu:consistance}, which  requires integrability of the
criterion, is always satisfied in the ballistic regime whereas it
might  fails  in  the  sub-ballistic  regime.  For  instance,  when  $
\sup_{\t' \in \Theta} |\dot \phi_{\t'}|$ is integrable with respect to
$\pi_\t$, point~\hyperref[assu:integrability]{$iii)$} of
Assumption~\ref{assu:consistance}  follows. This    occurs    in
Examples~\ref{ex:deuxpoints}        and~\ref{ex:beta}        presented
below. However, this point is not satisfied in Example~\ref{ex:Temkin}
as   suggested   by   point   (c)   of   Proposition~\ref{prop:temkin}
below. Nevertheless, we show the consistency of the MLE and prove that
the Fisher  information is infinite in this  framework suggesting that
the rate of convergence is faster than $\sqrt{n}$.

%%%%%%%%%%%%%%%%%%%%%%%%%%%%%%%%%%%%%%%%%%%%%%%%%
%%% Environment with finite and known support %%%
%%%%%%%%%%%%%%%%%%%%%%%%%%%%%%%%%%%%%%%%%%%%%%%%%

\begin{exam} \label{ex:deuxpoints}
Fix $a_1 < a_2 \in (0,1)$ and let
$\nu_p=p\delta_{a_1}+(1-p)\delta_{a_2}$, where  $\delta_a$ is the Dirac
mass located at value $a$. Here, the unknown parameter is the proportion
$p\in  \Theta  \subset  [0,1]$  (namely~$\theta=p$). We  suppose  that
$a_1$,      $a_2$      and       $\Theta$      are      such      that
Assumption~\ref{assu:rho} is satisfied.  
\end{exam}

This example  is easily generalized  to $\nu$ having $m\ge  2$ support
points namely $ \nu_\t=\sum_{i=1}^m p_ia_i$, where $a_1,\dots,a_m$ are
distinct,     fixed     and    known     in     $(0,1)$,    we     let
$p_m=1-\sum_{i=1}^{m-1}p_i$     and    the     parameter     is    now
$\theta=(p_1,\dots,p_{m-1})$.

In the framework of Example~\ref{ex:deuxpoints}, we have
\begin{equation} \label{equa:Phi2pts}
 \phi_{p}(x,y)      =\log    [    p    a_1^{x+1}(1-a_1)^y    +(1-p)
 a_2^{x+1}(1-a_2)^y ] ,
\end{equation}
 \begin{prop} \label{prop:2ptsNormal}
In  the framework  of  Example~\ref{ex:deuxpoints}, assuming  moreover
that  $\Theta   \subset  (0,1)$,  Assumptions  ~\ref{assu:consistance}
and~\ref{assu:nonsingular}  are satisfied,  and hence  the MLE  of the
parameter $p$ is consistent and asymptotically normal.  
\end{prop}

%%%%%%%%%%%%
%%% BETA %%%
%%%%%%%%%%%%

\begin{exam} \label{ex:beta}
We let $\nu_\t$ be  a Beta distribution with parameters $(\alpha,\beta)$,
namely
\[
\dd\nu_\t(a) = \frac 1 {\Beta(\alpha,\beta)} a^{\alpha
  -1}(1-a)^{\beta -1} \dd a, \quad
   \Beta(\alpha,\beta)   =   \int_0^1   t^{\alpha
  -1}(1-t)^{\beta -1}\dd t.
\]
Here,  the  unknown parameter  is  $\theta=(\alpha,\beta) \in  \Theta$
where $\Theta$ is a compact subset of 
\[
  \{ (\alpha,\beta) \in (0,+ \infty)^2 \, : \, \beta < \alpha \leq \beta +1\}.
\]
\end{exam}
%%%%%%%%%%%%%%%%%%%%%%%%%%
%%%commentaires exemple%%%
%%%%%%%%%%%%%%%%%%%%%%%%%%%
The  inequalities $\beta<\alpha$  and $\alpha  \leq \beta  +1$ ensures
that                                points~\hyperref[assu:trans]{$ii)$}
and~\hyperref[assu:sousbal]{$iii)$}  of  Assumption~\ref{assu:rho} are
satisfied.

In the framework of Example~\ref{ex:beta}, we have
\begin{equation}
  \phi_\t(x,y) = \log \frac{\Beta(x+1+\alpha,y+\beta)}{\Beta(\alpha,\beta)}
\end{equation}
\begin{prop} \label{prop:betaNormal}
In    the     framework    of    Example~\ref{ex:beta},    Assumptions
~\ref{assu:consistance} and~\ref{assu:normality} are satisfied, 
and hence  MLE of  the parameter $(\alpha,  \beta)$ is  consistent and
asymptotically normal. 
\end{prop}

%%%%%%%%%%%%%%%%%%%%%%%
%%% SUPPORT INCONNU %%%
%%%%%%%%%%%%%%%%%%%%%%%

\begin{exam}[Temkin model] \label{ex:Temkin}
  We let  $\nu_\t =  p \delta_{a} +  (1-p) \delta_{1-a}$, where  $p$ is
 fixed in $(0,1/2)$ and the  unknown parameter is $\theta=a \in \Theta$,
 where $\Theta$ is a compact subset of $(0,p)$.
\end{exam}
The inequalities $p<1/2$ and $a<p$ ensures that points~\hyperref[assu:trans]{$ii)$}
and~\hyperref[assu:sousbal]{$iii)$}  of  Assumption~\ref{assu:rho} are
satisfied.  

In this framework, we have 
\begin{equation} \label{equa:Phi_temkin}
Q_\t(u,v)= \binom{u+v}{v}\ee^{\phi_{\t}(u,v)} = p K_a(u,v) + (1-p) K_{1-a}(u,v),
\end{equation}
with $K_a(u,v)$ defined as
\begin{equation}
\label{equa:Ka}
K_a(u,v)=\binom{u+v}{v}a^{u+1}(1-a)^v.
\end{equation}
%%

%%%%%%%%%%%%%%%%%%%%%%%%%%%%%%%
%%% iii) N'EST PAS VERIFIEE %%%
%%%%%%%%%%%%%%%%%%%%%%%%%%%%%%%

\begin{prop} \label{prop:temkin}
In the framework of Example~\ref{ex:Temkin}, the following holds.

(a) For any $\alpha >0$,
\begin{equation} \label{equa:ergoTemkin}
\frac  1  n   \sum_{x=0}^{n-1} \ \sup_{\t  \in  \VV^\complement_\alpha}
\left[    \phi_\t(L^{n}_{x+1},L^n_x)   -   \phi_\ts(L^{n}_{x+1},L^n_x)
\right]  \xrightarrow[n   \to  \infty]{}  -   \infty,  \quad  \mbox{in
  $\PPs$-probability,} 
\end{equation}
where $\VV^\complement_\alpha$  is the complement of $\VV_\alpha$ defined as $
\VV_\alpha= \{ a \in \Theta \, : \, \KL(a^\star | a) \leq \alpha \}$, 
with $\KL(\cdot | \cdot)$ is the Kullback-Leibler distance on $(0,1) \times
(0,1)$ defined as
\[ \KL(q|q')  = q \log  \frac q {q'}  + (1-q) \log \frac  {1-q} {1-q'}
\geq 0.
\]
Therefore, the MLE of the parameter $a$ is consistent.

(b) The Fisher information is infinite, that is, for any $\t$,
\begin{equation}
\Sigma_{\t} = \Et \big[ ( \phi_{\t}')^2 \big] = + \infty.
\end{equation}
(c) For any $\t \neq
\ts$, 
\begin{equation} \label{equa:der_phi_inf}
  J=\Es \big[ |  \phi_{\t}'| \big] = + \infty.
\end{equation}
\end{prop}

% From (a), we deduce that the MLE of the parameter $a$ is
% consistent. Indeed, taking $\theta_0=\ts$ does not change the value of
% $ \htet$ and obviously $\lsb_n(\ts)=0, \forall n\in\N.$
% From \eqref{equa:ergoTemkin},
% $\displaystyle    \sup_{\t     \in    \VV_\alpha^\complement}    \frac
% 1n{\lsb_n(\theta)}\xrightarrow [n\to\infty]{}-\infty,$ 
% hence the consistency follows.  %%

%%%%%%%%%%%%%%
%%% PROOFS %%%
%%%%%%%%%%%%%%

\section{Proofs} \label{sect:proofs}

%%%%%%%%%%%%%%%%%%%%%%%
%%% LIMIT CRITERION %%%
%%%%%%%%%%%%%%%%%%%%%%%

\subsection{Proof of Proposition~\ref{prop:LimitEllTheta}} 
\label{sect:LimitEllTheta}
First, we establish the weak law of large numbers
\begin{equation} \label{equa:WLLN}
 \frac  1n  \lsb_n(\t) \xrightarrow[n  \to  \infty]{} \lsb(\t),  \quad
 \mbox{in $\PPs$-probability}. 
\end{equation}
Since the sequence $L_n^n$, $L_{n-1}^n$, $\dots$, $L^n_0$ has the same
distribution as the BPIRE $Z_0$, $Z_1$, $\dots$, $Z_n$ defined by
\eqref{eq:BPRE}, we have 
\begin{equation} \label{equa:ell_equi}
\lsb_n(\theta) \sim
 \sum_{k=0}^{n-1}\left[\phi_{\theta}(Z_k,Z_{k+1})
   -\phi_{\theta_0}(Z_k,Z_{k+1}) \right], 
\end{equation}
under $\PPs$, where $\sim$ means equality in distribution. 
\citeauthor{Comets_etal}           proved          that          under
point~\hyperref[assu:transience]{$ii)$}  of Assumption~\ref{assu:rho},
the  process  $(Z_n,Z_{n+1})_{n\in  \Z_+}$  is  a  positive  recurrent
homogeneous Markov chain which admits the unique invariant probability
measure $\tilde \pi_{\ts}$ defined by~\eqref{equa:pi_tilde}.  
Hence, according  to Theorem 4.2  in Chapter 4  from~\cite{Revuz}, for
any function $g : \Z_+^2\to 
\R^d$ such that $ \tilde \pi_{\ts}( \|g\|) <\infty$, the following ergodic
theorem holds 
\begin{equation} \label{equa:ergo}
  \lim_{n\to\infty}\frac 1 n \sum_{k=0}^{n-1} g(Z_k,Z_{k+1})
= \tilde \pi_{\ts}( g) ,
\end{equation}
$\PPs$-almost   surely    and   in   $    \mathbb{L}^1(\PPs)$.   Under
point~\hyperref[assu:integrab]{$iii)$}                               of
Assumption~\ref{assu:consistance},  we can  use~\eqref{equa:ergo} with
$g=\phi_\t  - \phi_{\tz}$,  and  combining with~\eqref{equa:ell_equi},
this yields~\eqref{equa:WLLN}. 

Now we turn to the local  uniform weak law of large numbers.This could
be verified  by the  same arguments  as in the  proof of  the standard
uniform law of  large numbers \cite[see Theorem 6.10  and its proof in
Appendix 6.A in][]{bierens} where  \eqref{equa:WLLN} plays the role of
the  weak law  of large  numbers  for a  random sample  in the  former
reference.

Indeed,        under        point~\hyperref[assu:cont]{$i)$}        of
Assumption~\ref{assu:consistance},  the map  $\t  \mapsto \phi_{\t}  -
\phi_{\tz}$          is         continuous,          and         under
point~\hyperref[assu:integrab]{$iii)$} of
Assumption~\ref{assu:consistance}, we have 
\[
\tilde  \pi_\ts  \Big(  \sup_{\t  \in  \Theta}  \left  |  \phi_{\t}  -
  \phi_{\tz} \right | \Big) < +\infty, 
\]
which implies that
\[
\tilde \pi_\ts \Big( \sup_{\t \in \Theta} \phi_{\t} - \phi_{\tz} \Big)
< +\infty \quad \mbox{and} \quad 
\tilde \pi_\ts \Big( \inf_{\t \in \Theta} \phi_{\t} - \phi_{\tz} \Big)
> -\infty. 
\]
Therefore, the proof of Theorem  6.10 in \cite{bierens} can be adapted
to our context and this implies~\eqref{equa:ULLN}. \hfill $\qed$

%%%%%%%%%%%%%%%%%%%%%
%% IDENTIFICATION %%%
%%%%%%%%%%%%%%%%%%%%%

\subsection{Proof of Proposition~\ref{prop:IdentificationEllTheta}} 
\label{sect:IdentificationEllTheta}

First       of       all,       note      that       under       under
point~\hyperref[assu:integrab]{$iii)$}                               of
Assumption~\ref{assu:consistance}, the limit $\lsb(\t)$ is finite for 
any value $\t \in \Theta$. 
%%Now,  we   start  by  proving  that  for   any  $  \theta\in\Theta,$
%%$\theta\neq \ts,$ we have $\ 
%%\lsb(\theta)< \lsb(\ts).$ 
From~\eqref{eq:defell}, we may write 
\begin{equation*}
  \lsb(\theta)-\lsb(\ts) = \tilde\pi_{\ts}(\phi_{\theta}-\phi_{\ts}).
\end{equation*}
Using~\eqref{equa:pi_tilde}    and   noting   that    $Q_{\t}(u,v)   =
\binom{u+v}{u} \exp[\phi_\t(u,v)]$ yields 
\[
 \lsb(\theta)-\lsb(\ts)=  \sum_{u \in \Z_+}
 \pi_{\ts}(u) \left[ \sum_{v \in \Z_+} \log 
 \left( 
 \frac{Q_\t(u,v)}{ Q_{\ts}(u,v)} 
 \right )
 Q_{\ts}(u,v)
 \right].
\]
Using  Jensen's inequality  with  respect  to the logarithm
function and the (conditional) distribution $Q_{\ts}(u,\cdot)$ yields 
\begin{equation}\label{equa:Jensen}
\lsb(\theta)-\lsb(\ts) \leq 
\sum_{u \in \Z_+} \pi_{\ts}(u) 
        \log                                                       
        \left[
        \sum_{v\in\Z_+}\frac{Q_\t(u,v)}{Q_{\ts}(u,v)       }Q_{\ts}(u,v)
        \right] =0.
\end{equation}
The equality in~\eqref{equa:Jensen}  occurs if and only if  for any $u
\in \Z_+$, we have $Q_\t(u,\cdot)=  Q_{\ts}(u,\cdot), $
which  is   equivalent  to  the  probability   measures  $\nu_\t$  and
$\nu_{\ts}$
having identical  moments.  Since their  supports are included  in the
bounded set 
$(0,1)$, these  probability measures are then  identical \cite[see for
instance][Chapter  II,  Paragraph 12,  Theorem  7]{Shir}.  Hence,  the
equality $\lsb(\t) =  \lsb(\ts)$ yields $\nu_{\theta}=\nu_{\ts}$ which
is equivalent to  $\t=\ts$ under point~\hyperref[assu:ident]{$ii)$} of
Assumption~\ref{assu:consistance}.  

In   other   words,  we   proved   that
$\lsb(\t)\le  \lsb(\ts)$ with  equality if  and only  if  $\t=\ts$. To
conclude  the proof  of Proposition~\ref{prop:IdentificationEllTheta},
it suffices to use that the function $\t 
\mapsto \lsb(\t)$ is continuous. \hfill $\qed$
%%%%%%%%%%%%%%%%%%%%%

%%%%%%%%%%%%%%%
%%% MOMENTS %%%
%%%%%%%%%%%%%%%

\subsection{Proof of Proposition~\ref{prop:pimoments}} \label{sect:pimoments}

Let  $\kappa_\t$  defined  by~\eqref{equa:kappa}  and  $\alpha$  be  a
positive number.  Let $\Lambda$ be the positive random variable such that 
\[
     1 - S =\ee^{-\Lambda},
\]
where $S$ is defined by~\eqref{equa:S}. Then, we have
\begin{equation} \label{equa:pi_lambda}
  \sum_{k=0}^{\infty}  k^{\alpha}   \pi_\t  (k)=\E^{\theta}  \left[  S
    \sum_{k=0}^{\infty} k^{\alpha} \ee^{-\Lambda k} \right]. 
\end{equation}
From the fact that for any integer $k$ and any positive $\lambda$
\[
\int_k^{k+1}x^{\alpha}  \ee^{-\lambda   x}\dd  x  \geq  \ee^{-\lambda}
k^{\alpha}\ee^{-\lambda k} 
\quad \mbox{and} \quad 
\int_k^{k+1}x^{\alpha}  \ee^{-\lambda  x}\dd
x\leq \ee^{\lambda} (k+1)^{\alpha} \ee^{-\lambda (k+1)}, 
\]
we deduce that 
\begin{equation} \label{equa:momentAlpha}
(1-S)\cdot  \frac {\Gamma(\alpha+1)}{ \Lambda^{1+\alpha}}
\leq 
\sum_{k=1}^{\infty}k^{\alpha} \ee^{-\Lambda k}
\leq 
\frac 1 {1-S} \cdot \frac {\Gamma(\alpha+1)}{ \Lambda^{1+\alpha}},
\end{equation}
where $\Gamma(z)=\int_0^{+\infty} x^{z-1} \ee^{-x} \dd x$. 
Using the fact that there exists a constant $C$ such that 
\[
  S \sum_{k=1}^{\infty}k^{\alpha} \ee^{-\Lambda k}  \1_{\{S>1/2\}} \leq C,
\]
that $\Lambda \geq S$ and~\eqref{equa:momentAlpha} yields 
\begin{equation} \label{equa:majo_pi_lambda}
  \Et \left[ S \sum_{k=0}^{\infty} k^{\alpha} \ee^{-\Lambda k} \right]
  \leq
   C + 2 \Gamma(\alpha+1) \Et \left[ S^{-\alpha} \right].
\end{equation}
\cite{Kesten} showed that there exists a positive constant $c_\t$ such
that
\begin{equation} \label{equa:Kesten}
\Pt(S^{-1}>x) \cdot x^{\kappa_\t} \to c_\t, \quad \mbox{when $x \to \infty$.}
\end{equation}
Combining~\eqref{equa:pi_lambda},~\eqref{equa:Kesten} 
and~\eqref{equa:majo_pi_lambda} implies point~\hyperref[prop:pimomentfini]{$i)$} of
Proposition~\ref{prop:pimoments}.

% Recall that
% %%
% \[
%   \Et  \left[ S^{-\alpha}  \right]  = \Et  \left[ \Big(  \sum_{k\geq0}
%     \prod_{i=1}^k \rho_i \Big)^{\alpha} \right]. 
% \]
% %%
% If $\alpha  \in (0,1]$, the  sub-additivity of the  function $s\mapsto
% s^{\alpha}$ combined  with the fact that  the sequence $(\rho_i)_{i\in
%   \N}$ is i.i.d. yields 
% %%
% \[
%   \Et \left[  \Big( \sum_{k\geq0} \prod_{i=1}^k  \rho_i \Big)^{\alpha}
%   \right]  \leq \Et  \left[ \sum_{k\geq0}  \Big(  \prod_{i=1}^k \rho_i
%     \Big)^{\alpha}  \right]= \sum_{k\geq0} \left(  \Et[ \rho_0^\alpha]
%   \right)^k. 
% \]
% %%
% If $\alpha >1,$ using Minkowski inequality, monotony and the fact that
% the sequence $(\rho_i)_{i\in \N}$ is i.i.d. yields 
% %%%%%%%%%%%%
% \[
%   \Et \left[  \Big( \sum_{k\geq0} \prod_{i=1}^k  \rho_i \Big)^{\alpha}
%   \right]^{1/\alpha} \leq 
%   \sum_{k\geq0} \Et  \left[ \Big( \prod_{i=1}^k  \rho_i \Big)^{\alpha}
%   \right]^{1/\alpha} 
%   = 
%   \sum_{k\geq0} \left( \Et[ \rho_0^\alpha]^{1/\alpha}\right)^k.
% \]
% %%
% Hence, for any positive  $\alpha$, the fact that $\Et[\rho_0^\alpha] <
% 1$                 combined                with~\eqref{equa:pi_lambda}
% and~\eqref{equa:majo_pi_lambda}                                 implies
% point~\hyperref[prop:pimomentfini]{$i)$}                             of
% Proposition~\ref{prop:pimoments}. 

%%%%%%%%%%%%%%%%
%%% 2ND STATEMENT %%%
%%%%%%%%%%%%%%%%%

Now,   we   turn   to   point~\hyperref[prop:pimomentfini]{$ii)$}   of
Proposition~\ref{prop:pimoments}. Using the  convexity of the function
$x \mapsto | \log(1-x) |$ on $(0,1)$, we obtain 
\[
  \frac{ \1_{\{S<1/2\}} } {2S\log 2 } \leq  \frac{ \1_{\{S<1/2\}} } {\Lambda}, 
\]
which combined with~\eqref{equa:momentAlpha} yields 
\[
  \Et \left[ S \sum_{k=0}^{\infty} k^{\alpha} \ee^{-\Lambda k} \right]
  \geq
   \frac{\Gamma(\alpha+1)}     {2     (2\log2)^{\alpha+1}    }     \Et
   \left[ S^{-\alpha} \1_{\{S<1/2\}} \right], 
\]
and finally
\begin{equation} \label{equa:mino_pi_lambda}
   \Et \left[ S \sum_{k=0}^{\infty} k^{\alpha} \ee^{-\Lambda k} \right]
   \geq
   \frac{\Gamma(\alpha+1)}   {2   (2\log2)^{\alpha+1}   }  \Big(   \Et
   \left[S^{-\alpha} \right] - 2^\alpha \Big). 
\end{equation}
Combining~\eqref{equa:pi_lambda},~\eqref{equa:Kesten}
and~\eqref{equa:mino_pi_lambda} implies point~\hyperref[prop:pimomentfini]{$ii)$} of
Proposition~\ref{prop:pimoments}.
% If  $\alpha \geq  1$,  the sur-additivity  of  the function  $s\mapsto
% s^{\alpha}$ combined  with the fact that  the sequence $(\rho_i)_{i\in
%   \N}$ is i.i.d. yields 
% %%
% \[
%   \Et \left[  \Big( \sum_{k\geq0} \prod_{i=1}^k \rho_i \Big)^{\alpha} \right] 
%   \geq 
%   \Et \left[  \sum_{k\geq0} \Big(  \prod_{i=1}^k \rho_i \Big)^{\alpha} \right]
%   = 
%   \sum_{k\geq0} \left( \Et[ \rho_0^\alpha] \right)^k.
% \]
% %%
% If  $\alpha \in  (0,1)$,  the  concavity of  the  function $s  \mapsto
% s^\alpha$  and the  fact  that the  sequence  $(\rho_i)_{i\in \N}$  is
% i.i.d. yields 
% %%
% \[
%   \Et \left[  \Big( \sum_{k\geq0} \prod_{i=1}^k \rho_i \Big)^{\alpha} \right]
%   \geq
%   \Et \left[  \Big( \sum_{k=0}^N \prod_{i=1}^k \rho_i \Big)^{\alpha} \right]
%   \geq
%   (N+1)^{\alpha-1} \sum_{k=0}^N \left( \Et[ \rho_0^\alpha] \right)^k,
% \]
% %%
% for any integer  $N$. Hence, for any positive  $\alpha$, the fact that
% $\Et[\rho_0^\alpha]   \geq  1$   combined  with~\eqref{equa:pi_lambda}
% and~\eqref{equa:mino_pi_lambda} implies
% point~\hyperref[prop:pimomentfini]{$ii)$} of
% Proposition~\ref{prop:pimoments}.
\hfill $\qed$

%%%%%%%%%%%%%%%%%
%%% EXAMPLE 1 %%%
%%%%%%%%%%%%%%%%%

\subsection{Proof of Proposition~\ref{prop:2ptsNormal}}
\citeauthor{FLM}       have       already       established       that
points~\hyperref[assu:cont]{$i)$}  and~\hyperref[assu:ident]{$ii)$} of
Assumption~\ref{assu:consistance}       as      well       as  
point~\hyperref[assu:diff]{$i)$}   of  Assumption~\ref{assu:normality}
are satisfied. From the latter  reference, we also know that the first
derivative $  \dot \phi_{p}$ as  well as the second  derivative $\ddot
\phi_{p}$  are uniformly  bounded when  $\Theta \in  (0,1)$,  and this
implies     that     point~\hyperref[assu:integrability]{$iii)$}    of
Assumption~\ref{assu:consistance} and
points~\hyperref[assu:smoothness]{$ii)$}
and~\hyperref[assu:uniform]{$iv)$}  of Assumption~\ref{assu:normality}
are       satisfied.      Points~\hyperref[assu:invertibility]{$iii)$}
and~\hyperref[assu:nonsingular]{$v)$}                                of
Assumption~\ref{assu:normality}   can  be   checked   exactly  as   in
\cite{FLM}. \hfill $\qed$

%%%%%%%%%%%%%%%%%
%%% EXAMPLE 2 %%%
%%%%%%%%%%%%%%%%%

\subsection{Proof of Proposition~\ref{prop:betaNormal}}
\citeauthor{FLM}       have       already       established       that
points~\hyperref[assu:cont]{$i)$}  and~\hyperref[assu:ident]{$ii)$} of
Assumption~\ref{assu:consistance}       as      well       as     
point~\hyperref[assu:diff]{$i)$}   of  Assumption~\ref{assu:normality}
are satisfied.  

From the latter reference, we  know that there exists a constant $A_1$
independent of $\theta$, such that for any $u$ and $v$ 
\begin{equation} \label{equa:majo_dot_phi_beta}
|\partial_\alpha \phi_\t (u,v)| \leq A_1 \log(1+v) \quad \mbox{and} \quad
|\partial_\beta \phi_\t (u,v)| \leq A_1\log(1+u).
\end{equation}
Define $\kappa_\t \in (0,1]$  as the unique positive number satisfying
$\Et [ \rho_0^{\kappa_\t} ]=1$, that is, 
\[
  \Gamma(\alpha-\kappa_\t)  \Gamma(\beta+\kappa_\t)  =  \Gamma(\alpha)
  \Gamma(\beta).  
\]
 Define  $\kappa=\min  \{  \kappa_\t  \,:\,  \t  \in  \Theta\}$.  From
 \eqref{equa:majo_dot_phi_beta},  there exists  $A_2>0$  and $A_3  >0$
 independent of $\theta$, such that for any $u$ and $v$ 
\begin{equation}  \label{equa:majo_dot_phi_beta2}
  |\partial_\alpha\phi_{\t}(u,v)| \leq A_2 v^{\kappa/2} 
  \quad \mbox{and} \quad
  |\partial_\beta\phi_{\t}(u,v)| \leq A_2 u^{\kappa/2}, \\
\end{equation}
and
\begin{equation}  \label{equa:majo_dot_phi_beta3}
  |\partial_\alpha\phi_{\t}(u,v)|^{4}\leq A_3v^{\kappa/2} 
  \quad \mbox{and} \quad
  |\partial_\beta\phi_{\t}(u,v)|^{4}\leq A_3 u^{\kappa/2}.
\end{equation}
Using  the  fact that  $\Et[  \rho_0^{\kappa/2}]<1$  for  any $\t  \in
\Theta$, Proposition~\ref{prop:pimoments}, the fact that 
\[
\sum_{k\in\Z_+}k^{\kappa/2}\pi_{\t}(k)       =       \sum_{u,v\in\Z_+}
u^{\kappa/2}  \tilde   \pi_\t(u,v)  =  \sum_{u,v\in\Z_+}  v^{\kappa/2}
\tilde \pi_\t(u,v), 
\]
\eqref{equa:majo_dot_phi_beta2}   and  \eqref{equa:majo_dot_phi_beta3}
yields     that     point~\hyperref[assu:integrability]{$iii)$}     of
Assumption~\ref{assu:consistance}    is   satisfied,   as    well   as
point~\hyperref[assu:smoothness]{$ii)$}                              of
Assumption~\ref{assu:normality} with $q=2$. 

Now,   we  turn   to   point~\hyperref[assu:invertibility]{$iii)$}  of
Assumption~\ref{assu:normality}. To exchange the order of derivation 
and summation, it is sufficient to prove that
\begin{equation} \label{equa:ConvNorm}
\sum_v \sup_{\t\in\Theta}\| \dot Q_\t (u,v) \| < \infty,
\end{equation}
for any integer $u$. 
Define $\theta'=(\alpha',\beta')$ with
\[
  \alpha' = \inf (\proj_1(\Theta)) \quad \mbox{and} \quad 
\beta' = \inf (\proj_2(\Theta)),
\]
where $\proj_i$, $i=1,2$ are the two projectors on the coordinates. 
Note that $\theta'$ does not necessarily belong to $\Theta$.
However, 
it still belongs to the sub-ballistic region. From \cite{FLM}, we know
that there exists a constant $A_4$ such that 
\[
  Q_{\t}(u,v) \leq A_4 Q_{\theta'}(u,v),
\]
for  any integers  $u$ and  $v$.  Define  $\kappa' \in  (0,1]$  as the
unique positive  number satisfying $\E^{\t'}  [ \rho_0^{\kappa'} ]=1$,
and recall that $\dot Q_\t(u,v) =
Q_\t(u,v) \dot \phi_\t (u,v)$. Hence, using the last inequality and 
the  fact  that  $\|  \dot  \phi_\t  \|  =  O(v^{\kappa'/2})$,  it  is
sufficient to prove that 
\begin{equation} \label{equa:MomentQ}
\sum_v v^{\kappa'/2} Q_{\theta'}(u,v) < \infty, \quad \mbox{for any integer $u$},
\end{equation}
to get \eqref{equa:ConvNorm}. We have
\[
  \sum_u \Big(\sum_v v^{\kappa'/2}  Q_{\theta'}(u,v) \Big) \pi_{\theta'}(u) =
  \sum_v v^{\kappa'/2} \pi_{\theta'}(v) < \infty,
\]
where  the  last  inequality   comes  from  the  fact  that  $\E^{\t'}
[\rho_0^{\kappa'/2}]<1$  and Proposition~\ref{prop:pimoments}.  Hence,
\eqref{equa:MomentQ} is 
satisfied for any integer  $u$ which proves that \eqref{equa:ConvNorm}
is satisfied.

The second order derivatives of $\phi_{\t}$ are given by 
 \begin{align*}
  \partial_{\alpha}^2  \phi_{\t}(x,y)   &  =  -\sum_{k=0}^x   \frac  1
  {(k+\alpha)^2} +\sum_{k=0}^{x+y} \frac 1 {(k+\alpha+\beta)^2} , \\
\partial_{\alpha}  \partial_{\beta} \phi_{\t}(x,y)  & =\sum_{k=0}^{x+y}
\frac 1 {(k+\alpha+\beta)^2} ,
\end{align*}
and  similar formulas  for  $\beta$ instead  of  $\alpha$.  Thus,  the
second derivative $\ddot \phi_{\t}$  is uniformly bounded on $\Theta$,
and   this   implies   that  point~\hyperref[assu:uniform]{$iv)$}   of
Assumption~\ref{assu:normality}                                      is
satisfied.   Point~\hyperref[assu:nonsingular]{$v)$}   of   Assumption
\ref{assu:consistance} can be checked exactly as in \cite{FLM}. \hfill
$\qed$

%%%%%%%%%%%%%%%%%
%%% EXAMPLE 3 %%%
%%%%%%%%%%%%%%%%%

\subsection{Proof of point (a) of Proposition~\ref{prop:temkin}}

Fix $\alpha >0$. To prove \eqref{equa:ergoTemkin}, we show 
\begin{enumerate}
\item[i)]
$\tilde   \pi_\ts  \big(  \big[\sup_{\t   \in  \VV^\complement_\alpha}
\phi_\t - \phi_\ts \big]^- \big)=+\infty$, 
\item[ii)]
$\tilde \pi_\ts \big( \big[\sup_{\t \in \VV^\complement_\alpha}\phi_\t
- \phi_\ts \big]^+ \big) < + \infty$. 
\end{enumerate}
Indeed, under  points i) and ii),  we can apply an  ergodic theorem to
$(Z_n)$ which yields 
\[
\frac  1  n   \sum_{x=0}^{n-1}  \sup_{\t  \in  \VV^\complement_\alpha}
\left[   \phi_\t(Z_{x},Z_{x+1})   -  \phi_\ts(Z_{x},Z_{x+1})   \right]
\xrightarrow[n  \to  \infty]{}  -  \infty,  \quad  \mbox{$\PPs$-almost
  surely,} 
\]
and then \eqref{equa:ergoTemkin}.  

We  note  that  $K_a(u,\cdot)$   defined  by  \eqref{equa:Ka}  is  the
distribution of a negative  binomial random variable $\NB(u+1,a)$ with
probability of  success $1-a$ and  number of failures $u+1$,  that is,
the  distribution  of  the  number  of  successes  in  a  sequence  of
independent Bernoulli trials until $u+1$ failures has occurred.  

We will make use several times of the fact that
$\NB(u+1,a)$ is the sum  of $(u+1)$ i.i.d.  geometric random variables
$G_1(a)$, $\dots$, $G_{u+1}(a)$ with parameter $1-a$, that is
$\mathrm{Prob}(G_1(a)=k)=(1-a)^k a$, whose mean is given by
$\mu=(1-a)/a>1$.  As a shortand of notation, we write $\mu^\star$ as the
ratio $(1-a^\star)/a^\star>1$. 

Define for any $\eps>0$ and any integer $u$, the sets 
\begin{align}
 \label{equa:A}
A(\eps,u) &= \Big\{  v \in \Z_+ \, : \,  \Big| \frac{v}{u+1} - \mu^\star
\Big| \leq \eps \Big\}, \\ 
\label{equa:B}
B(\eps,u)  &=  \Big\{  v  \in   \Z_+  \,  :  \,  \Big|  \frac{v}{u+1}  -
\frac{1}{\mu^\star} \Big| \leq \eps \Big\},\\ 
\label{equa:C}
C(\eps,u) &= \Z_+ \setminus \big( A(\eps) \cup B(\eps) \big).
\end{align}
We have
\[
\sum_{v \in A(\eps,u)} K_{a^\star}(u,v) =
\mathrm{Prob}\Big( |\overline{G}_{u+1}(a^\star) - \mu^\star | \leq \eps \big), 
\]
with
\[
\overline{G}_{u+1}(a^\star) = \frac 1 {u+1} \sum_{k=1}^{u+1} G_k(a^\star). 
\]
Using concentration inequalities, there exists a constant $c_\eps$ such that 
\begin{equation} \label{equa:probaA}
\sum_{v \in A(\eps,u)} K_{a^\star}(u,v)   \geq 1 - \ee^{c_\eps' (u+1)}.
\end{equation}

Similarly, there exists a constant $c_\eps'$ such that 
\begin{equation} \label{equa:probaB}
\sum_{v \in B(\eps,u)} K_{1-a^\star}(u,v) \geq 1 - \ee^{c_\eps'' (u+1)},
\end{equation}
and  as a consequence  of \eqref{equa:Phi_temkin},~\eqref{equa:probaA}
and~\eqref{equa:probaB}, there exists a constant $c_\eps''$ such that 
\begin{equation}
\label{equa:probaC}
\sum_{v \in C(\eps,u)} Q_{\ts}(u,v) \leq \ee^{-c''_\eps (u+1)}.
\end{equation}
Introduce the  quantity $\beta$ to be used later and defined as
\begin{equation}
\beta  = \max \Big\{\sup_{a  \in \VV_\alpha^\complement}  \Big| \log
\frac{1-a^\star}{1-a^{\phantom{\star}}} \Big|, \sup_{a \in \VV_\alpha^\complement} \Big|
\log \frac{a^\star}{a^{\phantom{\star}}}\Big| \Big\}. 
\end{equation}
For any $0<\eps< \mu^\star-1$ and for any $v$ in $A(\eps,u)$, we have $v
> u+1$ and as a consequence, for any $a \in \Theta$, 
\[
K_a(u,v) > K_{1-a}(u,v).
\]
Thus, we deduce that for any $v$ in $A(\eps,u)$,
\begin{align*}
(\phi_\t-\phi_\ts)(u,v)&=\log
\frac{pK_a(u,v)+(1-p)K_{1-a}(u,v)} {pK_{a^\star}(u,v)+(1-p)K_{1-a^\star}(u,v)}\\ 
&\leq
\log   \frac  1{p}   +   (u+1)  \log   \frac{a^{\phantom{\star}}}{a^\star}   +  v   \log
\frac{1-a^{\phantom{\star}}}{1-a^\star} \\ 
&\leq
\log  \frac  1{p} -  \frac{u+1}{a^\star}  \left[  \KL(a^\star  | a)  +
  \Big(\frac{v}{u+1}     -     \mu^\star     \Big)    a^\star     \log
  \frac{1-a^\star}{1-a^{\phantom{\star}}} \right]\\ 
&\leq 
\log \frac  1{p} - (u+1)  \Big( \frac{\alpha^{\phantom{\star}}}{a^\star} -  \eps \beta
\Big). 
\end{align*}
Similarly, we deduce that for any $0<\eps<1-1/\mu^\star$ and for any $v$ in
$B(\eps,u)$, 
\begin{align*}
(\phi_\t-\phi_\ts)(u,v)
&\leq
\log  \frac  1{1-p}  -  (u+1) \Big(  \frac{\alpha}{1-a^\star}  -  \eps
\beta\Big). 
\end{align*}
Hence, choosing 
\[
\eps=               \min\Big\{\mu^\star-1,              1-1/\mu^\star,
\frac{\alpha}{2\beta(1-a^\star)} \Big\} 
\] 
yields the existence of $u_0$ such that for any $u\geq u_0$ and any
$v$ in $A(\eps,u) \cup B(\eps,u)$
\begin{equation} \label{equa:logQ}
\Big(  \sup_{\t  \in \VV_\alpha^\complement}  (\phi_\t-\phi_\ts)(u,v)
\Big)^+=0 
\quad \mbox{and} \quad 
\Big( \sup_{\t \in \VV_\alpha^\complement} (\phi_\t - \phi_\ts)(u,v) \Big)^-\geq
\frac{\alpha}{3(1-a^\star)}(u+1). 
\end{equation}
Combining~\eqref{equa:probaC} and~\eqref{equa:logQ} immediatly yields
\[
\tilde \pi_\ts \big( \sup_{\t \in \VV_\alpha^\complement}\left[\phi_\t
  - \phi_\ts \right]^- \big) \geq \frac{\alpha}{3(1-a^\star)} 
\sum_{u \geq u_0} \pi_\ts(u) (u+1)(1-\ee^{-c''_\eps (u+1)}) = +\infty,
\]
where        the         last        equality        comes        from
Proposition~\ref{prop:pimoments}.    This   achieves   the  proof   of
point~i).  To prove  point  ii),  note that  there  exists a  positive
constant $c_1$ such that for any $u$ and any $v$
\[
\Big(  \sup_{\t  \in \VV_\alpha^\complement}  (\phi_\t-\phi_\ts)(u,v)
\Big)^+ \leq \Big| \sup_{\t \in \VV_\alpha^\complement} (\phi_\t - \phi_\ts)(u,v)
\Big| \leq c_1 (u+1+v).  
\]
Furthermore,   from   Cauchy-Schwarz   inequality,   the   fact   that
$K_{a^\star}(u,\cdot)$  (resp.  $K_{1-a^\star}(u,\cdot)$)  possesses a
second  moment  quadratic  with  $u$,  and~\eqref{equa:probaC},  there
exists two positive constants $c_2$ and $c_3$ such that 
\begin{align*}
\sum_{v \geq 0} Q_{\ts}(u,v) v \cdot \1_{C(\eps,u)} (v) 
&\leq
\Big(\sum_{v  \geq 0}  Q_{\ts}(u,v) v^2  \Big)^{1/2}\cdot \Big(\sum_{v
  \geq 0} Q_{\ts}(u,v) \1_{C(\eps,u)} (v)  \Big)^{1/2} \\
&\leq c_2 (u+1) \ee^{-c_3(u+1)},
\end{align*}
Therefore, there exists two positive constants $c_4$ and
$c_5$ such that 
\[
\tilde    \pi_\ts   \big(    \sup_{\t    \in   \VV_\alpha^\complement}
\left[\phi_\t - \phi_\ts \right]^+ \big) \leq c_4 + c_5 
\sum_{u \geq u_0} \pi_\ts(u) (u+1) \ee^{-c''(u+1)} < + \infty,
\]
which achieves the proof of point ii). 

Noting  that  $\htet$ does  not  depend on  the  choice  of $\t_0$  in
\eqref{eq:l},  we   can  take  $\theta_0=\ts$.    Obviously,  we  have
$\lsb_n(\ts)=0$,    for     all    integer    $n$,     whereas    from
\eqref{equa:ergoTemkin},   we  have   $\lsb_n(\t)/n$  which   goes  to
infinity, for any $\t$ outside a neighborhood of $\ts$. Hence, the
consistency follows. \hfill $\qed$ 

%%%%%%%%%%%%%%%%%%%%%%%%%%
%%% PROOF OF POINT (B) %%%
%%%%%%%%%%%%%%%%%%%%%%%%%%

\subsection{Proof of point (b) of Proposition~\ref{prop:temkin}}

We have,
 \begin{align} 
   \label{equa:phi_derivee}
   \phi_{\t}'(u,v) \cdot Q_\t(u,v)  &= p K_a(u,v) \left(\frac{u+1}{a} -
     \frac{v}{1-a}\right)    - (1-p) K_{1-a}(u,v)
   \left(\frac{u+1}{1-a} - \frac{v}{a}\right). 
\end{align}

Recall that 
\begin{equation*}
  \Sigma_{\t}=\Et  \big[   (  \phi_{\t}')^2  \big]   =  \sum_{u\geq0}
  \pi_{\t}(u) \sum_{v\geq0} Q_{\t}(u,v) [   \phi_{\t}'(u,v)]^2,
\end{equation*}
which can be rewritten using~\eqref{equa:phi_derivee} as
\begin{align*}
  \Sigma_{\t}     &=    \sum_{u\geq0}     \pi_{\t}(u)    \sum_{v\geq0}
  \frac{1}{Q_\t(u,v)} \left[ p K_a(u,v) 
   \left(\frac{u+1}{a}  -  \frac{v}{1-a}\right)  - (1-p)  K_{1-a}(u,v)
   \left(\frac{u+1}{1-a} - \frac{v}{a}\right) \right]^2.  
\end{align*}
Define
\[
v(u)= u+1, \quad 
V(u)=  \max\{  v \in  \Z_+  \, :  \,  v  \leq \mu  \cdot
(u+1)- (1-a)\sqrt{u+1}\},
\]
with  $\mu=(1-a)/a$. From the  fact that  $\mu>1$, there  exists $u_0$
such that for any 
$u\geq u_0$, we have  $v(u)<V(u)$. Furthermore, for any $u\geq u_0$
and any $v$ in $[v(u),V(u)]$, we have 
\[
K_a(u,v) \geq K_{1-a}(u,v), \quad
  \frac{u+1}{1-a} - \frac{v}{a} \leq 0 
\quad \mbox{and} \quad
\frac{u+1}{a} - \frac{v}{1-a} \geq \sqrt{u+1}.
\]
Thus,
\begin{align}
\nn   \Sigma_{\t}   &\geq   \sum_{u\geq  u_0}   \pi_{\t}(u)   \sum_{v=
  v(u)}^{V(u)} \frac{1}{Q_\t(u,v)}  \left[ p K_a(u,v)
    \Big(\frac{u+1}{a} - \frac{v}{1-a}\Big) \right]^2 \\ 
\label{equa:mino_I1}&\geq p^2 \sum_{u\geq u_0} (u+1)\pi_{\t}(u) \sum_{v=
    v(u)}^{V(u)} K_a(u,v).
\end{align} 
Recall that
\begin{align*}
\sum_{v= v(u)}^{V(u)} K_a (u,v) &= \mathrm{Prob} \Big( \NB(u+1,a) \in
 [v(u),V(u)] \Big) \\
&=\mathrm{Prob}\Big( \sqrt{u+1} \big( \overline{G}_{u+1} - \mu \big) \in
\big[\sqrt{u+1}(1-\mu), -(1-a) \big] \Big), 
\end{align*}
with
\[
\overline{G}_{u+1} = \frac 1 {u+1} \sum_{k=1}^{u+1} G_k(a), 
\]
where  $G_1(a)$, $\dots$,  $G_{u+1}(a)$ are  i.i.d.   geometric random
variables with mean $\mu$. 
From  the central limit  theorem applied  to the  sequence $(G_k(a))$,
there exists $u_1 \geq u_0$ such that for any $u \geq u_1$
\begin{equation}
  \label{equa:mino_Q1}
  \sum_{v= v(u)}^{V(u)} K_a(u,v) \geq \frac 1 2
  \mathrm{Prob}\Big( \NN(0,\sigma^2) \in \big[ -2 \mu, -(1-a) \big] 
  \Big) = C >0,
\end{equation}
where  $\sigma^2=(1-a)/a^2$   is  the  variance   of  $G_1(a)$  and
$\NN(0,\sigma^2)$  a  Gaussian random  variable  with  mean $0$  and
variance          $\sigma^2$.         Injecting~\eqref{equa:mino_Q1}
in~\eqref{equa:mino_I1} yields
\[
\Sigma_{\t} \geq C p^2 \sum_{u\geq u_1} (u+1) \pi_{\t}(u).
\]
From  Proposition~\ref{prop:pimoments}, $\pi_{\t}$  does not  possess a
finite first moment in the sub-ballistic regime and we 
deduce that $\Sigma_{\t}=+\infty$. \hfill $\qed$

%%%%%%%%%%%%%%%%%%%%%%%%%%
%%% PROOF OF POINT (C) %%%
%%%%%%%%%%%%%%%%%%%%%%%%%%

\subsection{Proof of point (c) of Proposition~\ref{prop:temkin}}

Assume that $\t \neq \ts$. Recall that
\begin{equation*}
  J=    \sum_{u\geq0}    \pi_{\ts}(u)    \sum_{v\geq0}    Q_{\ts}(u,v)
  |\phi_{\t}'(u,v) |, 
\end{equation*}
which can be rewritten using~\eqref{equa:phi_derivee} and~\eqref{equa:Q}
\begin{align*}
  J       &=       \sum_{u\geq0}       \pi_{\ts}(u)       \sum_{v\geq0}
  \frac{Q_\ts(u,v)}{Q_\t(u,v)} \left| p K_a(u,v) \left(\frac{u+1}{a} -
      \frac{v}{1-a}\right) -  (1-p) K_{1-a}(u,v) \left(\frac{u+1}{1-a}
      - \frac{v}{a}\right) \right|.  
\end{align*}
Using the set $A(\eps,u)$, we have 
\begin{equation}
J \geq p J_1 - (1-p) J_{2},
\end{equation}
with
\begin{align}
\label{equa:J1}
  J_1 &= \sum_{u\geq0} \pi_{\ts}(u) \sum_{v\in A(\eps,u)} 
   \frac{Q_\ts(u,v)}{Q_\t(u,v)}  K_a(u,v)
   \left|\frac{u+1}{a} - \frac{v}{1-a}\right| \\
\label{equa:J2}
  J_2 &=\sum_{u\geq0} \pi_{\ts}(u) \sum_{v\in A(\eps,u)} 
  \frac{Q_\ts(u,v)}{Q_\t(u,v)}  K_{1-a}(u,v)
   \left|\frac{u+1}{1-a} - \frac{v}{a}\right|.
\end{align}
Choose $\eps=\min\{|\mu-\mu^\star|/2,\mu^\star-1\}$. Then, 
for any $u$ and any $v \in A(\eps,u)$, we have
\[
\left|\frac{u+1}{a} -  \frac{v}{1-a}\right| \geq \frac{u+1}{1-a} \eps,
\quad K_a(u,v) \geq Q_\t(u,v),
\]
and as a consequence
\begin{align*}
J_1 & \geq  \frac{\eps}{1-a} \sum_{u\geq0} \pi_{\ts}(u)  (u+1) \sum_{v\in
  A(\eps,u)} Q_\ts(u,v) \\
&\geq \frac{p\eps}{1-a} \sum_{u\geq0} \pi_{\ts}(u)  (u+1) \sum_{v\in
  A(\eps,u)} K_{a^\star}(u,v).
\end{align*}
Using~\eqref{equa:probaA} and the fact that $\pi_\ts$  does not  possess a
finite first moment in the  sub-ballistic regime, we deduce that $J_1$
is infinite. 
On the other hand, we have for any $u$ and any $v \in A(\eps,u)$,
\[
\left|\frac{u+1}{1-a}  -  \frac{v}{a}\right|  \leq  (u+1)\Big(\frac  1
{1-a} + \frac{\eps+\mu^\star}{a}\Big),
\]
and,
\[
    \frac{K_{1-a}(u,v)}{Q_\t(u,v)} \leq \frac{K_{1-a}(u,v)}{pK_a(u,v)} 
\leq \frac 1 {p} \mu^{u+1-v}, 
\]
and as a consequence,
\begin{align*}
J_2 &\leq \frac{1}{p} \Big(\frac 1 {1-a} + \frac{\eps+\mu^\star}{a}\Big)
\sum_{u\geq0} \pi_{\ts}(u) (u+1) \cdot \gamma^{u+1} \sum_{v\in A(\eps,u)} Q_\ts(u,v),
\end{align*}
where
\[
\gamma= \mu^{-(\mu^\star-1-\eps)}<1.
\]
From  the fact that  $u \mapsto (u+1)
\gamma^{u+1}  \sum_{v\in  A(\eps,u)}  Q_\ts(u,v)  $  is  bounded,  and
therefore  integrable  against  $\pi_\ts$,  we deduce  that  $J_2$  is
finite.  This  achieves the proof  of~\eqref{equa:der_phi_inf}. \hfill
$\qed$

%%%%%%%%%%%%%%%%%%%
%%% SIMULATIONS %%%
%%%%%%%%%%%%%%%%%%%

\section{Numerical performance}
\label{sect:simus}
%{Simulation experiment: comparing estimation procedures}

In this section, we explore the numerical performance of our 
estimation procedure in  the frameworks of Example~\ref{ex:deuxpoints}
and the Temkin model. We compare our performance with the performance of
the   estimator  proposed   by   \cite{AdEn}. An explicit description of the form
of~\citeauthor{AdEn}'s  estimator  in   the  particular  case  of  the
one-dimensional nearest  neighbour path is provided in  Section 5.1 of
\cite{Comets_etal}. Therefore, one can  estimate  $\ts$ by  the  solution of  an
appropriate system of equations, as illustrated below.  

\textbf{Example~\ref{ex:deuxpoints} (continued).}  In this case the
parameter $\theta$ equals $p$ and we have 
\[
  v   =     \Es[\omega_0] = p^\star a_1 + (1-p^\star)a_2.
\]
Hence, among the visited sites, the proportion of those from which the
first  move is  to the  right gives  an estimator  for $p^\star  a_1 +
(1-p^\star)a_2$.  Using this observation, we can estimate $p^\star$. 

\textbf{Example~\ref{ex:Temkin} (continued).} In this case the
parameter $\theta$ equals $a$ and we have 
\[
  v   =     \Es[\omega_0] = p a^\star + (1-p)(1-a^\star).
\]
Hence, among the visited sites, the proportion of those from which the
first  move is  to  the right  gives  an estimator  for  $p a^\star  +
(1-p)(1-a^\star)$.  Using this observation, we can estimate $a^\star$.

\subsection{Experiments}
\label{sect:experiments}
We   now   present   the   simulation  experiment   corresponding   to
Example~\ref{ex:deuxpoints}   and  Example~\ref{ex:Temkin}   where  we
include a comparison with \citeauthor{AdEn}'s procedure.
     
For  each  of  the  two  simulations, we  \textit{a  priori}  fix  a
parameter  value $\ts$ as  given in  Table~\ref{tabl:theta_values} and
repeat 1,000 times the  procedure described below. 
\renewcommand{\arraystretch}{1.5}% Wider
\begin{table}[h]
  \centering
  \begin{tabular}{|c|c|c|}
    \hline
Simulation & Fixed parameter & Estimated parameter \\
    \hline
Example~\ref{ex:deuxpoints}   &   $(a_1,a2)=(0.4,0.7)$,  $\kappa=0.9$&
$p^\star \approx 0.548$\\ 
    \hline
Example~\ref{ex:Temkin} & $\kappa=0.9$, $p \approx 0.41$ & $a^\star=0.4$ \\
    \hline
  \end{tabular}
  \caption{Parameter values for each experiment.}
  \label{tabl:theta_values}
\end{table}

Then, we generate a random environment according to $\nu_\ts$ on the
set of  sites $\{-10^3,$ $\dots,  10^3\}$.  In fact,  we do not  use the
environment values for  all the $10^3$ negative sites,  since only few
of these sites are visited  by the walk.  However the computation cost
is very low comparing to the rest of the estimation procedure, and the 
symmetry is convenient for programming purpose.  Then, we run a random
walk in this environment and stop it successively at the hitting times
$T_n$ defined by  \eqref{equa:HittingTime}, with $n \in \{10^2  k \, :
\, 1\le k \le 10 \}$.  For each 
stop, we estimate $\ts$ according to our procedure and 
\citeauthor{AdEn}'s  one.  The  likelihood optimization  procedure was
performed  as a combination  of golden  section search  and successive
parabolic interpolation.  

The parameter is chosen such that the RWRE is transient to the right and
sub-ballistic. Note that the length of  the random walk is not $n$ but
rather  $T_n$.  The  fluctuations of  $T_n$  depend in  nature on  the
parameter  $\kappa$.  Under  mild additional  assumptions,  \cite{KKS}
proved that if $\kappa < 1$, then $n^{-1/\kappa} \cdot T_n$ has a
non-degenerate limit distribution, a stable law with index $\kappa$. 

In the simulations, the quantity $T_n$ varies considerably. To avoid too
long computations, when  $T_n$ is too large, we  fixed a threshold for
the number of steps for the walk at $t_{\max}= 500 n^{1/\kappa} \approx
10^6$.  When  the  threshold  is  reached,  we  did  not  compute  our
estimator. 
This case happened for $4.4\%$ (when $n=100$) and for $41.9\%$ (when
$n=1000$)  of the simulation  in Example~\ref{ex:deuxpoints},  and for
$0.3\%$ (when $n=100$) and for $4.9\%$ (when 
$n=1000$) of the simulation in Example~\ref{ex:Temkin}.

\begin{figure}[H]
  \centering
  \includegraphics[height=12.5cm,width=10cm,angle=-90]{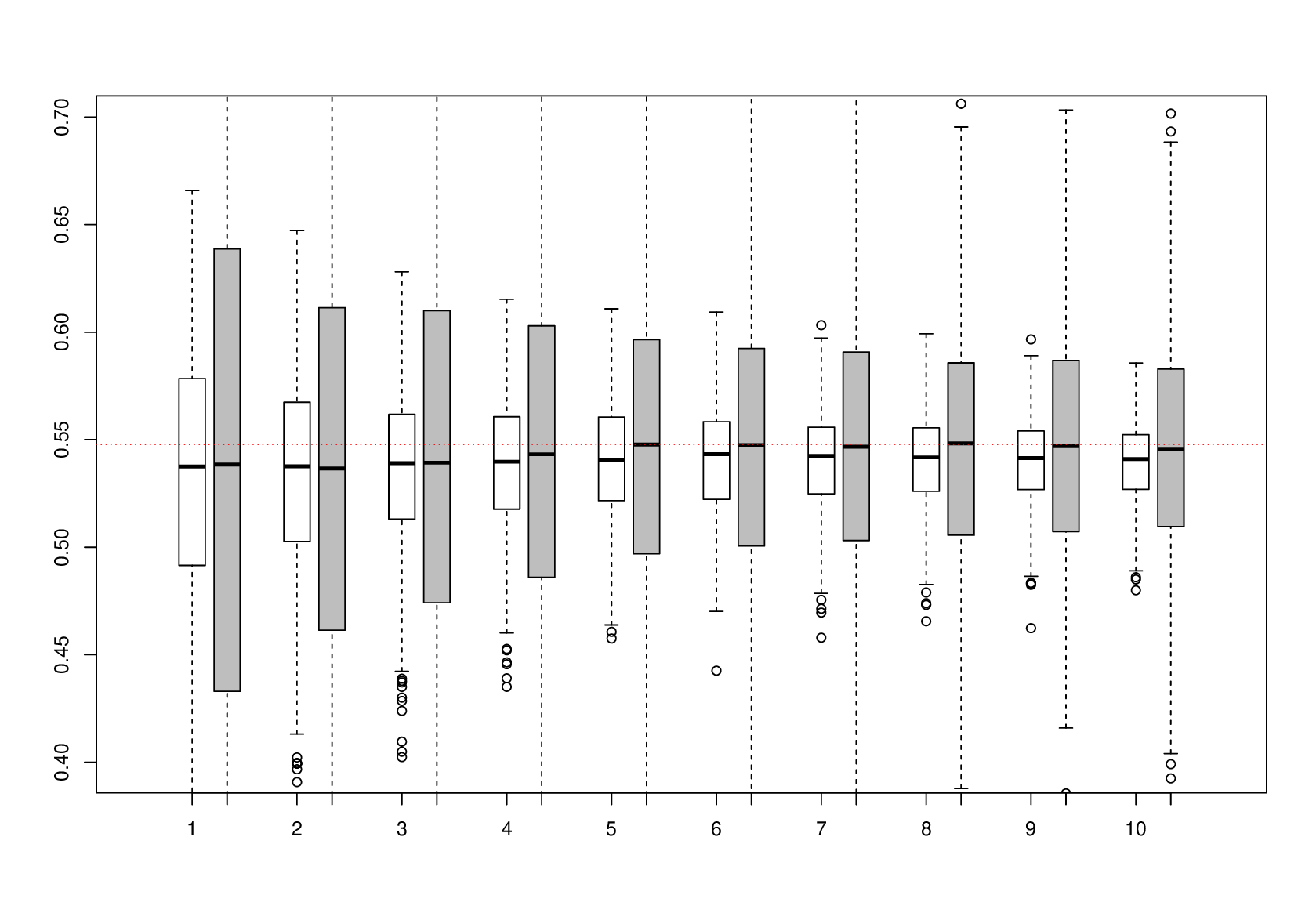}
  \\ 
  \caption{Boxplots   of   our   estimator   (left  and   white)   and
    \citeauthor{AdEn}'s estimator (right and grey) obtained from 1,000
    iterations and for values $n$ ranging
    in $  \{10^2 k ; 1\le k  \le 10 \}$ ($x$-axis  indicates the value
    $k$).  The panel displays estimation of $p^\star$ in
    Example~\ref{ex:deuxpoints}.  The true value  is indicated by horizontal lines.}
  \label{figu:boxplots_deuxpoints}
\end{figure}

Figure~\ref{figu:boxplots_deuxpoints}   shows  the  boxplots   of  our
estimator  and  \citeauthor{AdEn}'s   estimator  obtained  from  1,000
iterations of the procedures in Example~\ref{ex:deuxpoints}. First, we
shall  notify that  in  order  to simplify  the  visualisation of  the
results,    we   removed    in   the    boxplots    corresponding   to
Example~\ref{ex:deuxpoints} about 
1.5\% of outliers values (outside 1.5 times the interquartile range above
the   upper  quartile  and   below  the   lower  quartile)   from  our
estimator.  We observe that the accuracies of the
procedures increase with the value  of $n$.  We also note that whereas
\citeauthor{AdEn}'s seems unbiased our procedure seems to be slightly 
biaised. However, our procedure  exhibits a much smaller variance than
\citeauthor{AdEn}'s one.  One explanation for the worse performance of 
\citeauthor{AdEn}'s estimator  comparing to our procedure  is the fact
that only a few part of the trajectory is used in the estimation. 

\begin{figure}[H]
  \centering
  \includegraphics[height=12.5cm,width=10cm,angle=-90]{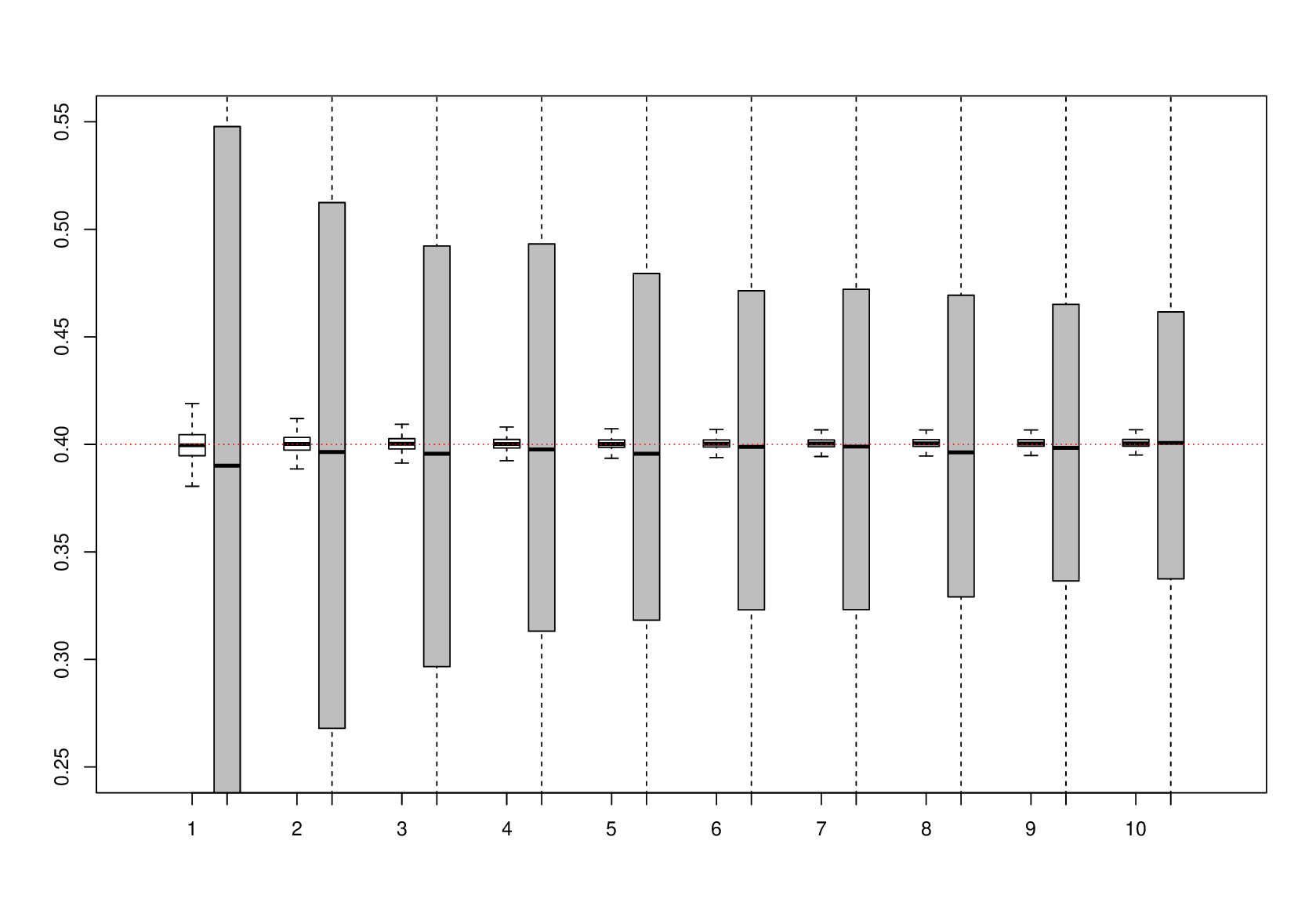} \\
  \caption{Boxplots   of   our   estimator   (left  and   white)   and
    \citeauthor{AdEn}'s estimator (right and grey) obtained from 1,000
    iterations and for values $n$ ranging
    in $  \{10^2 k ; 1\le k  \le 10 \}$ ($x$-axis  indicates the value
    $k$).  The panel displays estimation of $a^\star$ in
    Example~\ref{ex:Temkin}.  The true value  is indicated by horizontal lines.}
  \label{figu:boxplots_temkin}
\end{figure}

Figure~\ref{figu:boxplots_temkin} shows  the  boxplots   of  our
estimator  and  \citeauthor{AdEn}'s   estimator  obtained  from  1,000
iterations  of the procedures  in Example~\ref{ex:Temkin}. First, we
shall  notify that  in  order  to simplify  the  visualisation of  the
results,    we   removed    in   the    boxplots    corresponding   to
Example~\ref{ex:deuxpoints} about 
15\% of outliers values (outside 1.5 times the interquartile range above
the   upper  quartile  and   below  the   lower  quartile)   from  our
estimator.  We first
observe that the accuracies of the 
procedures increase  with the  value of $n$.   We also note  that both
procedures  seem  unbiased. However,  our  procedure  exhibits a  much
smaller variance than \citeauthor{AdEn}'s one, but also a much smaller
one than when  we were not estimating the  support. This suggests that
the  rate of  convergence when  estimating the  support in  the Temkin
model is faster than the square root of $n$.

%%%%%%
%%%%%%
\bibliographystyle{chicago}
\bibliography{MAMA_2013_11_21}

\end{document}